\documentclass[10pt,a4paper]{amsart}
\usepackage[latin1]{inputenc}
\usepackage{amsmath}
\usepackage{amsfonts}
\usepackage{amssymb}
\usepackage{graphicx}
\usepackage{url}
\usepackage{tikz}

\usepackage{listings}

\usepackage{amsthm}
\theoremstyle{plain}
\newtheorem{thm}{Theorem}[section]

\newtheorem{cj}{Conjecture}[section]

\def\R{\mathbb{R}}

\def\K{\mathbb{K}}
\def\P{\mathbb{P}}

\begin{document}

\title{On the Banach-Mazur Distance between the Cube and the Crosspolytope}
\author{Fei Xue}
\address{Institut für Mathematik, Technische Universit\"at Berlin,
  Sekr. Ma 4-1, Straße des 17 Juni 136, D-10623 Berlin, Germany}
\email{xue@math.tu-berlin.de}
\thanks{The research was supported by a PhD scholarship of the Berlin
  Mathematical School}

\begin{abstract}
	% In convex geometry we concern on the distance between different convex bodies. The Banach-Mazur distance indicates the distance between two symmetric convex bodies. However for given two convex bodies it is usually difficult to confirm the exact value of Banach-Mazur distance.
        In this note we study the Banach-Mazur distance between the
        $n$-dimensional cube and the crosspolytope. Previous work
        shows that the distance has order $\sqrt{n}$, and here we will
        prove some explicit bounds improving on former
        results. Even in dimension 3 the exact distance is not known,
        and based on computational results it is conjectured to 
        be $\frac{9}{5}$. Here we will also present computerbased potential optimal results in dimension $4$ to $8$.
\end{abstract}

\maketitle

\section{Introduction}

We call $K\subset \R^n$ an $n$-dimensional convex body, if $K$ is compact, and for any $x,y\in K$ and $\lambda\in [0,1]$, $\lambda x+(1-\lambda)y\in K$. The set of all $n$-dimensional convex bodies is denoted as $\K^n$. A convex polytope $P$ is defined as the convex hull of finitely many points
$$P=conv\{u_1,\cdots,u_k\}$$
and the set of all $n$-dimensional convex polytopes is denoted as $\P^n$.

The Hausdorff distance between two convex bodies $K$ and $L$ is defined as:
$$d_H(K,L)=\max\{\sup_{x\in K}\inf_{y\in L}d(x,y),\sup_{y\in L}\inf_{x\in K}d(x,y)\}$$
where $d(x,y)$ is the usual Euclidean distance. Equivalently,
$$d_H(K,L)=\inf\{\epsilon>0:X\subset Y_{\epsilon},Y\subset X_{\epsilon}\}$$
where $X_{\epsilon}=\cup_{x\in X}\{z:d(z,x)\leq\epsilon\}$.

For a real number $p\geq 1$, the $p$-norm of $x\in \R^n$ is defined by
$$||x||_p=(|x_1|^p+|x_2|^p+\cdots+|x_n|^p)^\frac{1}{p}.$$
The maximum norm is the limit of the $p$-norm for $p\rightarrow\infty$. It is equivalent to
$$||x||_\infty=\max{|x_1|,|x_2|,\cdots,|x_n|}.$$ 

Denote 
$$C_n=\{x\in\R^n:||x||_\infty\leq 1\}=[-1,1]^n$$
the $n$-dimensional unit cube, and 
$$C_n^\star=\{x\in\R^n:||x||_1\leq 1\}$$
the $n$-dimensional unit crosspolytope. Denote 
$$B_n=\{x\in\R^n:||x||_2\leq 1\}$$
the $n$-dimensional unit ball. For examples, the Hausdorff distance between $C_n$ and $C_n^\star$ is $\frac{n-1}{\sqrt{n}}$, and the Hausdorff distance between $C_n$ and $B_n$ is $\sqrt{n}-1$.

The Banach-Mazur distance between two symmetric convex bodies $K$ and $L$ is defined as:
$$d_{BM}(K,L)=\min\{r>0:K\subset gL\subset rK,g\in GL(n,\R)\}$$
where $GL(n,\R)$ is the group of linear transformations. It can be deduced that
$$d_{BM}(K_1,K_3)\leq d_{BM}(K_1,K_2)d_{BM}(K_2,K_3).$$
Due to this reason, the Banach-Mazur distance is also written as 
$$\log\min\{r:K\subset gL\subset rK,g\in GL(n,\R)\}.$$
In this paper we only concern about the distance, so we keep the definition to be the former one.

There are some results on the Banach-Mazur distance to some special convex bodies. John's Theorem on the maximal ellipsoid contained in a convex body gives the estimate:

\begin{thm}[John's Theorem\cite{Jo48}]
	The Banach-Mazur distance between an $n$-dimensional convex body $K$ and $n$-dimensional ball is at most $\sqrt{n}$.
\end{thm}

As a corollary, for any two convex bodies $K$ and $L$, 
$$d_{BM}(K,L)\leq d_{BM}(K,B_n)d_{BM}(B_n,L)\leq n.$$
As a matter of fact, the diameter of $(\K^n,d_{BM})$ is still unknown, but E.Gluskin \cite{Gl81} proved that the diameter is bounded below $c n$ for some universal $c>0$. 

For symmetric reasons, one can easily prove that:

\begin{thm}[\cite{Ve09}]
	The Banach-Mazur distance between $B_n$ and $C_n$ is $\sqrt{n}$. The Banach-Mazur distance between $B_n$ and $C_n^\star$ is $\sqrt{n}$.
\end{thm}

There are also some results on the Banach-Mazur distance from any convex body to the cube \cite{BS88,Gi95}.

We are interested in the Banach-Mazur distance between $C_n$ and $C_n^\star$ \cite{Web2}. There are results in \cite{TJ89,Ve09} showing that the distance has order $\sqrt{n}$:

\begin{thm}[\cite{TJ89,Ve09}]
	There exists constant $c,C>0$ such that
	$$c\sqrt{n}\leq d_{BM}(C_n,C_n^\star)\leq C\sqrt{n}.$$
\end{thm}

To be exact, for the upper bound one can get 
$$C=\frac{1}{\sqrt[4]{2}-1}=5.2852\cdots.$$
For the lower bound, the constant $C$ is not explicitly stated in \cite{TJ89}.

In this paper we discuss the upper bound and the lower bound of this distance. Our main results are:

\begin{thm}
	(1) There is a maximum absolute constant $\alpha$ independent of dimension $n$, such that for any $x\in\R^n$ with $||x||_2=1$,
	$$\frac{1}{2^n}\sum_{v_i\in\{-1,1\}^n}|<x,v_i>|\geq \alpha.$$\\
	(2) $\alpha>\frac{1}{1.71453\cdots}$.
\end{thm}

\begin{thm} With $\alpha$ from above we have 
	$$\alpha\sqrt{n}\leq d_{BM}(C_n,C_n^\star)\leq (\sqrt{2}+1)\sqrt{n}.$$
\end{thm}

\section{Some results by computer}

To find the Banach-Mazur distance between the cube and the crosspolytope, one need to find the optimal $g\in GL(n,\R)$ and minimum $r>0$ for 
$$\frac{1}{r}C_n\subset gC_n^\star\subset C_n.$$
Assume that $g$ is the linear transformation $T=(x_{ij})_{n\times n}$, then the crosspolytope
$$gC_n^\star=conv\{\pm(x_{i1},\cdots,x_{in}):i=1,\cdots,n\},$$ 
and $gC_n^\star\subset C_n$ implies that $|x_{ij}|\leq 1$ for $i,j=1,\cdots,n$. The left part $\frac{1}{r}C_n\subset gC_n^\star$ with miminum $r$ implies that the vertices of the cube $\frac{1}{r}C_n$ is contained in the crosspolytope $gC_n^\star$, which is 
$$\max_{v_i\in\{-1,1\}^n} ||T^{-1}.v_i||_{1}=r$$
where $T=(x_{ij})_{n\times n}$. Therefore the Banach-Mazur distance is
$$d_{BM}(C_n,C_n^\star)=\min_{T} \max_{v_i\in\{-1,1\}^n} ||T^{-1}.v_i||_{1}$$
where $T=(x_{ij})_{n\times n}$ with $|x_{ij}|\leq 1$.

In principle this problem can be solved by a computer programme like Maple and Mathematica. We can use the code here on Wolfram Mathematica:
\begin{lstlisting}[language=Mathematica,frame=shadowbox]
dim = 3;
T = Array[Subscript[TT, ##] &, {dim, dim}];
B1 = IdentityMatrix[dim];
B1 = Join[-B1, B1];
Binf = Tuples[{-1, 1}, dim];
NMinimize[
Join[{Max[Table[Norm[Inverse[T].Binf[[j]], 1], 
{j, Length[Binf]}]], Det[T] != 0}, 
Table[Norm[T.B1[[i]], Infinity] <= 1, {i, Length[B1]}]], 
Flatten[T]]
\end{lstlisting}
which is offered by \cite{Web1}, where we can change 3 to any dimension we need. Since the computer only gives the numerical results, we made some adjustment to make them to be the probably optimal ones.

In dimension $3$ the distance is $\frac{9}{5}$ and the crosspolytope is:
\begin{equation}
\left(
\begin{matrix}
1 & 1 & -1/3\\
-1/3 & 1 & 1\\
1 & -1/3 & 1
\end{matrix}
\right).
\end{equation}

In dimension $4$ the distance is $2$ and the crosspolytope is
\begin{equation}
\left(
\begin{matrix}
1 & 1 & 1 & -1\\
-1 & 1 & 1 & 1\\
1 & -1 & 1 & 1\\
1 & 1 & -1 & 1
\end{matrix}
\right).
\end{equation}

In dimension $5$ the distance is $2.32871$ and the crosspolytope is
\begin{equation}
\left(
\begin{matrix}
0.792559 & 1 & 0.0387439 & -1 & -0.704555\\
1 & 0.792092 & 0.999411 & 0.855944 & 1\\
-1 & -0.0773263 & 1 & -1 & 0.888962\\
0.925403 & -1 & 1 & -0.115724 &-0.822648\\
1 & -0.79255 & -0.999989 & -0.856439 & 1
\end{matrix}
\right).
\end{equation}
It seems to be highly irregular.

In dimension $6$ the distance is $2.4488$ and the crosspolytope is
\begin{equation}
\left(
\begin{matrix}
1 & 1 & 1 & 1 & 1 & 1\\
-1 & x & 1 & y & -1 & 1\\
-1 & 1 & x & 1 & y & -1\\
-1 & -1 & 1 & x & 1 & y\\
-1 & y & -1 & 1 & x & 1\\
-1 & 1 & y & -1 & 1 & x
\end{matrix}
\right)
\end{equation}
where $x=0.324842$, $y=-0.434446$.

In dimension $7$ the distance is $2.6$ and the crosspolytope is
\begin{equation}
\left(
\begin{matrix}
1 & 1 & 1 & 1 & 1 & 1 & 1\\
1 & 0 & 1 & -1 & -1 & 1 & -1\\
1 & 1 & 0 & 1 & -1 & -1 & -1\\
1 & -1 & 1 & 0 & 1 & -1 & -1\\
1 & -1 & -1 & 1 & 0 & 1 & -1\\
1 & 1 & -1 & -1 & 1 & 0 & -1\\
1 & -1 & -1 & -1 & -1 & -1 & 1
\end{matrix}
\right).
\end{equation}

In dimension $8$ the distance is $2.5$, smaller than in dimension $7$, and the crosspolytope derives from one Hadamard matrix:
\begin{equation}
\left(
\begin{matrix}
1&1&1&1&1&1&1&1\\
-1&-1&-1&1&-1&1&1&1\\
-1&1&-1&-1&1&-1&1&1\\
-1&1&1&-1&-1&1&-1&1\\
-1&1&1&1&-1&-1&1&-1\\
-1&-1&1&1&1&-1&-1&1\\
-1&1&-1&1&1&1&-1&-1\\
-1&-1&1&-1&1&1&1&-1
\end{matrix}
\right).
\end{equation}

\section{Upper bound}

Recall that the Banach-Mazur distance between the cube and the crosspolytope is 
$$d_{BM}(C_n,C_n^\star)=\min_{T} \max_{v_i\in\{-1,1\}^n} ||T^{-1}.v_i||_{1}$$
where $T=(x_{ij})_{n\times n}$ with $|x_{ij}|\leq 1$. By giving a special $T$ one can get an upper bound of the distance.

\subsection{Hadamard matrix}

A Hadamard matrix is a square matrix whose entries are either $+1$ or $-1$ and whose rows are mutually orthogonal. A Hadamard matrix has maximal determinant among matrices with entries of absolute value less than or equal to $1$.

Sylvester \cite{Sy67} provided one way to construct the Hadamard matrix. Let 
$$H_1=\left(1\right)$$
$$H_2=
\left(
\begin{matrix}
1 & 1\\
1 & -1
\end{matrix}
\right)$$
and
$$H_{2^k}=
\left(
\begin{matrix}
H_{2^{k-1}} & H_{2^{k-1}}\\
H_{2^{k-1}} & -H_{2^{k-1}}
\end{matrix}
\right)$$
for $k\geq 2$, then $H_{2^k}$ are all Hadamard matrices.

The Hadamard conjecture proposes that a Hadamard matrix of order $4k$ exists for every positive integer $k$. Sylvester's construction yields Hadamard matrices of order $2^k$. A generalization of Sylvester's construction proves that if $H_n$ and $H_m$ are Hadamard matrices of orders $n$ and $m$ respectively, then there exists Hadamard matrix of order $nm$. So far the Hadamard conjecture is still open.

\subsection{Upper bound}

We are going to prove that:

\begin{thm}
	$$d_{BM}(C_n,C_n^\star)\leq (\sqrt{2}+1)\sqrt{n}.$$
\end{thm}

In dimension $n=2^k$, there exists a Hadamard matrix $H_{2^k}$. Choose the matrix $T_{2^k}=H_{2^k}$, then $T_{2^k}^{-1}=\frac{1}{n}T_{2^k}^t$ where $T_{2^k}^t$ is still a Hadamard matrix with row vectors $T_1,\cdots,T_n$. So
\begin{eqnarray*}
	&&\max_{v_i\in\{-1,1\}^n} ||T_{2^k}^{-1}.v_i||_{1}\\
	&&=\max(|<T_1,v_i>|+\cdots+|<T_n,v_i>|)\\
	&&\leq \max \sqrt{n(<T_1,v_i>^2+\cdots+<T_n,v_i>^2)}\\
	&&=\max \sqrt{n\cdot\frac{1}{n}\cdot||v_i||_2^2}\\
	&&=\sqrt{n}.
\end{eqnarray*}

Assume that in dimension $t\leq 2^k$ the upper bound is not bigger than $(\sqrt2+1)\sqrt{t}$ with crosspolytope $T_{t}$. Then in dimension $n=2^k+t$ where $t\leq 2^k$, let 
$$T_{2^k+t}=
\left(
\begin{matrix}
T_{2^k} & 0 \\
0 & T_{t}
\end{matrix}
\right).$$ 
The distance is therefore
\begin{eqnarray*}
	&&\max_{v_i\in\{-1,1\}^n} ||T_{2^k+t}^{-1}.v_i||_{1}\\
	&&=\max_{v_i\in\{-1,1\}^{2^k}} ||T_{2^k}^{-1}.v_i||_{1}+\max_{v_i\in\{-1,1\}^t} ||T_{t}^{-1}.v_i||_{1}\\
	&&\leq \sqrt{2^k}+(\sqrt2+1)\sqrt{t}\\
	&&\leq (\sqrt2+1)\sqrt{2^k+t}\\
	&&=(\sqrt2+1)\sqrt{n}.
\end{eqnarray*}

\subsection{Extension}

The Hadamard conjecture asked for the existence of Hadamard matrix in dimension $n=4k$. When the Hadamard matrix exists in dimension $n=4k$, denoted as $H_{4k}$, choose the vertices of the crosspolytope to be the vertors of Hadamard matrix $H_{4k}$, and the distance will be $\sqrt{n}$.

When $n=4k+j$, let the vertices of the crosspolytope to be the vectors of 
$$\left(
\begin{matrix}
I_j & 0 \\
0 & H_{4k}
\end{matrix}
\right).$$
Then the distance is $\sqrt{4k}+j$. Therefore the upper bound will be $\sqrt{n}+3$ for all $n$.

\section{Lower bound}

Recall that $\alpha$ is an absolute constant introduced in Theorem 1.4. In this section we are going to prove that:

\begin{thm}
	$$d_{BM}(C_n,C_n^\star)\geq \alpha\sqrt{n}.$$
\end{thm}

The Banach-Mazur distance of the cube and the octahedron is to find the minimum value of
$$\max_{v_i\in\{-1,1\}^n}||T^{-1}.v_i||_{1}.$$

Without loss of generality, consider only $\det(T)>0$. Write $T^{-1}=\det(T^{-1})^{1/n}N$, where $N\in SL(n,\R)$, the group of special linear transformations. Let the row vectors of $N$ be $N_j$, i.e. $N=\big(N_j\big)_{n\times 1}$, then we have
$$||N.v_i||_1=|<N_1,v_i>|+\cdots+|<N_n,v_i>|.$$
Also, since $\det(N)=1$, by the definition of determinant we have:
$$\prod_{j=1}^{n}||N_j||_2\geq 1$$
and by the arithmetic geometric inequality
$$\sum_{j=1}^{n}||N_j||_2\geq n({\prod_{j=1}^{n}||N_j||_2})^{1/n}\geq n.$$

Recall that $\alpha$ is an absolute constant independent of dimension $n$, such that for any $x\in\R^n$ with $||x||_2=1$,
$$\frac{1}{2^n}\sum_{v_i\in\{-1,1\}^n}|<x,v_i>|\geq \alpha.$$
We will discuss the existence and the value of $\alpha$ in the next section. Since the left hand side is linear, for any $x\in\R^n$, we have
$$\frac{1}{2^n}\sum_{v_i\in\{-1,1\}^n}|<x,v_i>|\geq \alpha||x||_2.$$

Based on this result, we can infer that:
\begin{eqnarray*}
	&&\max_{v_i\in\{-1,1\}^n}||T^{-1}.v_i||_1\\
	&&=\det(T^{-1})^{1/n}\max_{v_i\in\{-1,1\}^n}||N.v_i||_1\\
	&&=\det(T^{-1})^{1/n}\max_{v_i\in\{-1,1\}^n}\sum_{j=1}^{n}|<N_j,v_i>|\\
	&&\geq \det(T^{-1})^{1/n}\frac{1}{2^n}\sum_{v_i\in\{-1,1\}^n}\sum_{j=1}^{n}|<N_j,v_i>|\\
	&&=\det(T^{-1})^{1/n}\frac{1}{2^n}\sum_{j=1}^{n}\sum_{v_i\in\{-1,1\}^n}|<N_j,v_i>|\\
	&&\geq \alpha\det(T^{-1})^{1/n}\sum_{j=1}^{n}||N_j||_2\\
	&&\geq \alpha\det(T^{-1})^{1/n}n(\prod_{j=1}^{n}||N_j||_2)^{1/n}\\
	&&\geq \alpha\det(T^{-1})^{1/n}n\\
	&&\geq \alpha\sqrt{n}.
\end{eqnarray*}
The last inequality comes from: since $|x_{ij}|\leq 1$, we have $\det(T)\leq n^{n/2}$.

\section{An average problem}

We are looking for the maximal absolute constant $\alpha$ such that
$$\frac{1}{2^n}\sum_{v_i\in\{-1,1\}^n}|<x,v_i>|\geq \alpha||x||_2$$
holds for all $x\in\R^n$ and all dimension $n$.

\begin{cj}
	$$\frac{1}{2^n}\sum_{v_i\in\{-1,1\}^n}|<x,v_i>|\geq
        \frac{1}{\sqrt{2}}||x||_2,$$
i.e., $\alpha=1/\sqrt{2}$.
\end{cj}

In this section, we will prove that:

\begin{thm}
	$$\frac{1}{2^n}\sum_{v_i\in\{-1,1\}^n}|<x,v_i>|> \frac{1}{1.71453}||x||_2.$$
\end{thm}

\subsection{A special class of polytope}

A convex polytope may be defined as an intersection of a finite number of half-spaces. Which is to say, for any convex polytope $P$, there exists vectors $u_j$ ($1\leq j\leq k$) such that
$$P=\{x\in\R^n:<x,u_j>\leq 1;1\leq j\leq k\}.$$
For the same reason, for any symmetric convex polytope $C$, there exists vectors $u_j$ ($1\leq j\leq k$) such that
$$P=\{x\in\R^n:|<x,u_j>|\leq 1;1\leq j\leq k\}.$$

Consider the set
$$K=\{x\in\R^n:\sum_{j=1}^{k}|<x,u_j>|\leq 1\}$$
where $u_j$ ($1\leq j\leq k$) are non-zero vectors such that $K$ is
bounded. As the intersection of $2^k$ halfspaces $K$ is a convex polytope. 

% For any two vectors $x,y\in\R^n$ and vector $u\in\R^n$, we have 
% $$\frac{|<x,u>|+|<y,u>|}{2}\geq\frac{|<x,u>+<y,u>|}{2}=|<\frac{x+y}{2},u>|.$$
% Therefore when $x,y\in K$, we have $\frac{x+y}{2}\in K$, meaning that $K$ is convex.

\subsection{Proof in dimension $n=2,3,4$}
For general dimension $n$, the problem is equivalent to find the maximal value of $||x||_2$ in the convex polytope
$$F_n(x)=\frac{1}{2^n}\sum_{v_i\in\{-1,1\}^n}|<x,v_i>|\leq 1.$$
The maximal value is the value on some special vertices of this convex polytope. Moreover, if $x$ is a vertex of this convex polytope, then it is the intersection of at least $n$ facets.  

We assert that there are $n-1$ linearly independent
$v_i\in\{-1,1\}^n$ such that $<x,v>=0$. 

For any $x,y\in\R^n$ and very small $\epsilon>0$,
$$F_n(x+\epsilon y)+F_n(x-\epsilon y)-2F_n(x)=2\sum_{v_i\in\{-1,1\}^n\atop <x,v_i>=0}|<\epsilon y,v_i>|.$$
Notice that if $<y,v_i>=0$ for all $v_i$ such that $<x,v_i>=0$, and if $F_n(x)=F_n(x+\epsilon y)=1$, then we have $F_n(x-y)=1$, which means that $x$ is not a vertex of the polytope.

For a point $x$ with $F_n(x)=1$, if there are at most $n-2$ linearly independent $v_i\in\{-1,1\}^n$ such that $<x,v>=0$, then there exists $y$ not linear to $x$, such that $<y,v_i>=0$ whenever $<x,v>=0$, meaning that
$$F_n(x+\epsilon y)+F_n(x-\epsilon y)-2F_n(x)=0$$
for any $\epsilon>0$. Now we choose $\epsilon$ to be arbitrary small, and let $y=y'+y''$, where $x+\epsilon y'$ is on some facet of the polytope containing $x$, while $y''$ is paralleled with $x$. Since $<y'',v_i>=0$ whenever $<x,v>=0$, we have
$$F_n(x+\epsilon y')+F_n(x-\epsilon y')=2F_n(x),$$
thus $x-\epsilon y'$ is also on some facet of the polytope. Therefore $x$ is not a vertex of the polytope.

With this observation, we can find out the vertices of the convex polytope $F_n(x)\leq 1$.

In dimension $2$, $F_2(x)=\frac{|x_1+x_2|+|x_1-x_2|}{2}\leq 1$ is the cube $C_2$, and the maximal value of $||x||_2$ is $\sqrt2$.

In dimension $3$, without loss of generality, assume that $x=(x_1,x_2,x_3)$, where $x_1\geq x_2\geq x_3\geq 0$. When $x_2+x_3\geq x_1$, we have
\begin{eqnarray*}
	F_3(x)&&=\frac{1}{4}(|x_1+x_2+x_3|+|x_1+x_2-x_3|+|x_2+x_3-x_1|+|x_3+x_1-x_2|)\\
	&&=\frac{x_1+x_2+x_3}{2}=1.
\end{eqnarray*}
When $x_1\geq x_2+x_3$, we have
\begin{eqnarray*}
	F_3(x)&&=\frac{1}{4}(|x_1+x_2+x_3|+|x_1+x_2-x_3|+|x_2+x_3-x_1|+|x_3+x_1-x_2|)\\
	&&=x_1=1.
\end{eqnarray*}

So the convex polytope is
$$conv\{(\pm1,\pm1,0),(\pm1,0,\pm1),(0,\pm1,\pm1)\}$$
and the maximum value of $||x||_2$ is $\sqrt2$.

In dimension $4$, consider the vertex $x=(x_1,x_2,x_3,x_4)$. If there are three linearly independent $v_i\in\{-1,1\}^4$ such that $<x,v>=0$, denoted as $v_1,v_2,v_3$, then:

(1)if $v_i$ and $v_j$ has 1 or 3 coordinate(s) in common, for example $v_1=(1,1,1,1)$ and $v_2=(1,1,1,-1)$, then $x_4=0$, and $(x_1,x_2,x_3)$ is a vertex of the polytope $F_3(x)\leq 1$.

(2)if all pairs of $v_i$ and $v_j$ have 2 coordinate in common, then $x$ has the form $(\pm t,\pm t,\pm t,\pm t)$, with (by calculation) $t=\frac{2}{3}$.

Therefore the convex polytope contains only one more series of vertices: $$(\pm\frac{2}{3},\pm\frac{2}{3},\pm\frac{2}{3},\pm\frac{2}{3}).$$ 
The maximal value of $||x||_2$ is still $\sqrt2$.

\subsection{Some classes of vertices}

The maximal value $||x||_2$ in the convex polytope
$$F_n(x)=\frac{1}{2^n}\sum_{v_i\in\{-1,1\}^n}|<x,v_i>|\leq 1$$
is the value on some special vertices. Moreover, if $x$ is a vertex of this convex polytope, then there are $n-1$ linearly independent $v_i\in\{-1,1\}^n$ such that $<x,v>=0$.

As a bad news, for the convex polytope in dimension $n$, there are much more vertices than the vertices in dimension $n-1$.

When $n=2k$ is even, the point $$\frac{(1,1,\cdots,1)}{F_n((1,1,\cdots,1))}$$
is a vertex of this convex polytope. Moreover, 
$$\frac{(1,\cdots,1,2j-1)}{F_n((1,\cdots,1,2j-1))}$$ 
($j=1,2,\cdots,k-1$) are also vertices of this convex polytope, because the vectors $(\pm 1,\cdots,\pm 1,-1)$ with $k+j-1$ times $1$ and $k-j$ times $-1$ are $(n-1)$-dimensional.

When $n=2k+1$ is odd, the points 
$$\frac{(1,\cdots,1,2j-2)}{F_n((1,\cdots,1,2j-2))}$$
($j=1,2,\cdots,k-1$) are vertices of the convex polytope, because the vectors $(\pm 1,\cdots,\pm 1,-1)$ with $k+j-1$ times $1$ and $k-j+1$ times $-1$ are $(n-1)$-dimensional.

There are also other kinds of vertices and we cannot enumerate all of them.

\subsection{A bound for the value}

We may assume that 
$$\frac{||x||_2}{F_n(x)}\leq \alpha_n$$ 
where
$$F_n(x)=\frac{1}{2^n}\sum_{v_i\in\{-1,1\}^n}|<x,v_i>|.$$
Then for any $x=(x_1,\cdots,x_n,x_{n+1})$, without loss of generality, assume that $x_1\geq x_2\geq \cdots\geq x_{n+1}\geq 0$. Let $y_{(1)}=(x_1,\cdots,x_n+x_{n+1})$, $y_{(2)}=(x_1,\cdots,x_n-x_{n+1})$. By definition we have
$$F_{n+1}(x)=\frac{F_n(y_{(1)})+F_n(y_{(2)})}{2}$$
and
$$||x||_2^2=\frac{||y_{(1)}||_2^2+||y_{(2)}||_2^2}{2}.$$
Therefore
\begin{eqnarray*}
	\frac{||x||_2}{F_{n+1}(x)}&&=\frac{\sqrt{2}\sqrt{||y_{(1)}||_2^2+||y_{(2)}||_2^2}}{F_n(y_{(1)})+F_n(y_{(2)})}\\
	&&\leq \alpha_n\frac{\sqrt{2}\sqrt{||y_{(1)}||_2^2+||y_{(2)}||_2^2}}{||y_{(1)}||_2+||y_{(2)}||_2}
\end{eqnarray*}
Recall that $x_1\geq x_2\geq \cdots\geq x_{n+1}\geq 0$, we have
\begin{eqnarray*}
	\frac{||y_{(1)}||_2^2}{||y_{(1)}||_2^2}&&=1+\frac{4x_nx_{n+1}}{x_1^2+\cdots+x_{n-1}^2+(x_n-x_{n+1})^2}\\
	&&\leq 1+\frac{4x_n^2}{(n-1)x_n^2}\\
	&&= 1+\frac{4}{n-1}
\end{eqnarray*}
Therefore by monotonicity we have
\begin{eqnarray*}
	\alpha_{n+1}&&=\frac{||x||_2}{F(x)}\leq \alpha_n\frac{\sqrt{2}\sqrt{2+\frac{4}{n-1}}}{1+\sqrt{1+\frac{4}{n-1}}}\\
	&&\leq \alpha_n(1+\frac{1}{2(n-1)^2}).
\end{eqnarray*}

Since we already know that $\alpha_4=\sqrt{2}$, by induction we have:
\begin{eqnarray*}
	\alpha_n &&\leq \sqrt{2}\prod_{j=4}^{n-1}(1+\frac{1}{2(j-1)^2})\\
	&&<\sqrt{2}\prod_{j=4}^{\infty}(1+\frac{1}{2(j-1)^2})\\
	&&\approx \sqrt{2}*1.21236\\
	&&\approx 1.71453.
\end{eqnarray*}

So we get:

\begin{thm}
	$$\frac{1}{2^n}\sum_{v_i\in\{-1,1\}^n}|<x,v_i>|> \frac{1}{1.71453}||x||_2.$$
\end{thm}

To get a better value, for example, if we can prove $\alpha_9=\sqrt{2}$ and use this induction, we can get $\alpha_n<1.50765\cdots$ If we can prove $\alpha_{k}=\sqrt{2}$, we can get
$$\alpha_n<\sqrt{2}\prod_{j=k}^{\infty}(1+\frac{1}{2(j-1)^2}).$$

So far we still believe that $\alpha_n=\sqrt{2}$ for all $n$.

\section{Extension of this average problem}

\subsection{Integration}

To find the maximal absolute constant $\alpha$ for
$$\frac{1}{2^n}\sum_{v_i\in\{-1,1\}^n}|<x,v_i>|\geq \alpha||x||_2,$$
we can also find the minimal value of $F_n(x)$ on the hypersphere $||x||_2=1$.

On the hypersphere $||x||_2=1$, the average value of $|<x,v_i>|$ is the fraction of the volume of $n$-dimensional half-hypersphere and the volume of $(n-1)$-dimensional hypersphere. Therefore one can get the average value of $F_n(x)$.

On the hypersphere $||x||_2=1$, the maximum value of $F(x)$ is naturally $1$. Therefore 
$$F_n(x+dx)\geq F_n(x)-F_n(dx)\geq F_n(x)-|dx|$$
and by integration we have
$$F_n(y)\geq F_n(x)-arc<x,y>.$$
If we have a series of points that are dense enough on the hypersphere, like $\epsilon$-net, we will get a lower bound depending on $\epsilon$ and the minimum value over the point set.

However, when $n\rightarrow\infty$, the set of the vertices of the cube $\frac{1}{\sqrt{n}}C_n$ is quite discrete on the hypersphere.

\subsection{Local Optimization}

On the hypersphere $||x||_2=1$, $F_n(x)$ is a continuous function, and $\alpha_n$ is the minimum value of $F_n(x)$. If there are not $n-1$ independent vertices in $\{-1,1\}^n$ vertical to $x$, we can always choose a tangent vector $dx$ such that 
$$F_n(x+dx)+F_n(x-dx)-2F_n(x)=0,$$
therefore at point $x$ it is not a local minimal value. So the local minimal value appears only when there are $n-1$ independent vertices in $\{-1,1\}^n$ vertical to $x$.

It shows the same observation as the aspect of convex polytope. 

%\subsection{Global Optimization}

%Since we can construct other kinds of deduction like $x\rightarrow y_1,y_2$, we can get other bounds even more close to $\sqrt2$. That's why we believe $\sqrt2$ is optimal.

%This problem is similar with Reinhardt conjecture: many local minimums in any dimensions, different values for different local minimums, paths between any two points, partial derivations with direction exist everywhere.

%\subsection{Torus}

%Assume that $X_i$ crosses $X$ with $x_i$ times. When we paste all $X_i$ together, what is the average crossing number with $X$?

%Too bad torus is just a model.

\section{Conclusion}

In this paper we give a new upper bound $(\sqrt{2}+1)\sqrt{n}$ and a conjectured better upper bound $\sqrt{n}+3$, and a new lower bound $\frac{\sqrt{n}}{1.71453}$ and a conjectured better lower bound $\frac{\sqrt{n}}{\sqrt{2}}$.


\begin{thebibliography}{}
	
	\bibitem{BS88}
	J.Bourgain, S.J.Szarek.
	\newblock The Banach-Mazur distance to the cube and the Dvoretzky-Rogers factorization.
	\newblock {\em Israel J. Math.}, 62:169, 1988.
	
	\bibitem{Gi95}
	A.A.Giannopoulos.
	\newblock A Note on the Banach-Mazur Distance to the Cube.
	\newblock {\em Geometric Aspects of Functional Analysis.} Operator Theory Advances and Applications, vol 77: 67-73, 1995.
	
	\bibitem{Gl81}
	E.D.Gluskin.
	\newblock The diameter of the Minkowski compactum is roughly equal to n (Russian).
	\newblock {\em Funktsional. Anal. i Prilozhen.}, 15(1):72?73, 1981.
	
	\bibitem{GrL87}
	P.M.Gruber, C.G.Lekkerkerker. 
	\newblock Geometry of Numbers.
	\newblock Amsterdam, Netherlands: North-Holland, 1987.
	
	\bibitem{Jo48}
	F.John.
	\newblock Extremum problems with inequalities as subsidiary conditions.
	\newblock {\em Studies and Essays Presented to R. Courant on his 60th Birthday}, Interscience Publishers, Inc., New York. January 8, 187-204, 1948.
	
	\bibitem{Sy67}
	J.J.Sylvester.
	\newblock Thoughts on inverse orthogonal matrices, simultaneous sign successions, and tessellated pavements in two or more colours, with applications to Newton's rule, ornamental tile-work, and the theory of numbers. 
	\newblock {\em Philosophical Magazine}, 34:461?475, 1867
	
	\bibitem{TJ89}
	N.Tomczak-Jaegermann.
	\newblock Banach-Mazur distances and finite-dimensional operator ideals.
	\newblock {\em Pitman monographs and surveys in pure and applied Mathematics.}, 38, 1989
	
	\bibitem{Ve09}
	R.Vershynin.
	\newblock Lectures in geometric functional analysis, 2009.
	\newblock \url{http://www-personal.umich.edu/~romanv/papers/GFA-book/GFA-book.pdf}.
	
	\bibitem{Web1}
	\url{https://gist.github.com/anonymous/8f2bbb5ddb9a463a592ae6e6bfe363d6}
	
	\bibitem{Web2}
	\url{http://mathoverflow.net/questions/237567/banach-mazur-distance-between-the-cube-and-the-octahedron}
	
	
\end{thebibliography}
\end{document}